\documentclass[12pt,a4paper]{article}
\usepackage[latin1]{inputenc}
\usepackage[T1]{fontenc}
\usepackage{amsmath,amssymb,amsthm,color}
\usepackage{graphicx}
\usepackage{amsopn}
\usepackage{color} 

\DeclareMathOperator*{\limm}{l.i.m}
\title{Mellin Analysis and Its Distance Concept   Applications to Sampling Theory}
\author{Carlo Bardaro  \\
 \small Department of Mathematics and Computer Sciences
 \small University of Perugia\\
 \small carlo.bardaro@unipg.it\\
 Paul L. Butzer\\
 \small Lehrstuhl A fuer Mathematik, RWTH Aachen,
 \small D-52056 Aachen, Germany,\\
 \small butzer@rwth-aachen.de\\
 Ilaria Mantellini\\
 \small Department of Mathematics and Computer Sciences,
 \small University of Perugia,\\
 \small mantell@dmi.unipg.it\\
 Gerhard Schmeisser\\
 \small Department of Mathematics, FAU Erlangen-Nuremberg
 \small D-91058 Erlangen, Germany\\
 \small schmeisser@mi.uni-erlangen.de}
 
\begin{document}
\maketitle

\noindent
{\small {\bf Abstract:} In this paper a notion of functional ``distance'' in the Mellin transform setting is introduced and a general representation formula
 is obtained for it. Also, a determination of the distance is given in terms of Lipschitz classes and Mellin-Sobolev spaces. Finally applications to 
approximate versions of certain basic relations valid for Mellin band-limited functions are studied in details.}
\section{Introduction}
A Mellin version of the Paley-Wiener theorem of Fourier analysis was introduced in \cite{BBMS}, using both complex and real approaches. Moreover, the 
structure of the set of Mellin band-limited functions (i.e. functions with compactly supported Mellin transform) was studied. It turns out that a Mellin 
band-limited function cannot at the same time be Fourier band-limited, and it is extendable as an analytic function to the Riemann surface of the (complex) 
logarithm. This makes the theory of Mellin band-limited functions very different from the Fourier band-limited ones since one has to extend the notion of 
the Bernstein spaces in a suitable way, involving Riemann surfaces (Mellin-Bernstein spaces). In the classical frame, Fourier band-limitedness is a very 
fundamental assumption in order to obtain certain basic formulae such as the Shannon sampling theorem, the Mellin reproducing kernel formula, the Boas 
differentiation formula, the Bernstein inequality, quadrature formulae and so on. When a function $f$ is no longer (Fourier) band-limited, certain 
approximate versions of the above formulae are available with a remainder which needs to be estimated in a suitable way. This was done in \cite{BSS1}, 
\cite{BSS2}, \cite{BSS3} in terms of an appropriate notion of ``distance'' of $f$ from the involved Bernstein space.  In the Mellin transform setting an 
exponential version of the Shannon sampling theorem for Mellin band-limited functions was first introduced in a formal way in \cite{OSP}, \cite{BP} in 
order to study problems arising in optical physics. A precise mathematical version of the exponential sampling formula, also in the approximate sense, was given in \cite{BJ3}, \cite{BJ2}, employing a rigorous Mellin transform analysis, as developed in 
\cite{BJ0}, \cite{BJ1} (see also \cite{BBM2}). Furthermore, a Mellin version of the reproducing kernel formula, both for Mellin band-limited funtions and 
in an approximate sense, was proved in \cite{BBM0}. Therefore it is quite natural to study estimates of the error in the approximate versions of the 
exponential sampling theorem, the reproducing kernel formula, the Bernstein inequality and the Mellin-Boas differentiation formula using a new notion of 
``Mellin distance'' of a function $f$ from a Mellin-Bernstein space. In the present paper, we introduce a notion of distance in the Mellin frame, 
 and we prove certain basic representation theorems for it (Sec.~3). In Sec.~4 we give precise evaluations of the Mellin distance in some fundamental 
function spaces such as Lipschitz classes and Mellin-Sobolev spaces. In Sec.~5 we describe some important applications to the approximate exponential 
sampling thoerem, the Mellin reproducing kernel theorem and the Boas differentiation formula in the Mellin setting, employing Mellin derivatives. Moreover, 
the theory developed here enables one to obtain an interesting approximate version of the Bernstein inequality with an estimation of the remainder. 
The present approach may also be employed in order to study other basic relations valid for Mellin band-limited functions.

\section{Notations and preliminary results}

Let $C(\mathbb{R}^+)$\,and $C(\{c\} \times i\mathbb{R})$ be the spaces of all uniformly continuous and bounded functions defined on $\mathbb{R}^+$ and on 
the line $\{c\} \times i\mathbb{R}, c \in \mathbb{R},$ respectively, endowed with the usual sup-norm $\|\cdot\|_\infty,$ and  let $C_0(\mathbb{R}^+)$  be the 
subspace of 
$C(\mathbb{R}^+)$ of functions $f$ satisfying $\lim_{x \rightarrow 0^+} f(x) = \lim_{x \rightarrow +\infty}f(x) = 0.$
For $1\leq p < +\infty,$ let $L^p= L^p(\mathbb{R}^+)$~ be the space of all the Lebesgue measurable and $p$-integrable complex-valued functions defined on 
$\mathbb{R}^+$ endowed with the usual norm $\|\cdot\|_p.$ Analogous notations hold for functions 
defined on $\mathbb{R}.$

For $p=1$ and $c \in \mathbb{R},$ let us consider the space
$$X^1_c = \{ f: \mathbb{R}^+\rightarrow \mathbb{C}: f(\cdot) (\cdot)^{c-1}\in L^1(\mathbb{R}^+) \}$$
endowed with the norm
$$ \| f\|_{X^1_c} := \|f(\cdot) (\cdot)^{c-1} \|_1 = \int_0^\infty |f(u)|u^{c-1} du.$$

More generally, let $X^p_c$  denote the space of all  functions $f: \mathbb{R}^+\rightarrow \mathbb{C}$ such that $f(\cdot) (\cdot)^{c-1/p}\in 
L^p(\mathbb{R}^+)$ with $1<p< \infty.$ In an equivalent form, $X^p_c$ is the space of all functions $f$ such that $(\cdot)^c f(\cdot) \in 
L^p_\mu(\mathbb{R}^+),$ where $L^p_\mu= L^p_\mu(\mathbb{R}^+)$ denotes the Lebesgue space with respect to the (invariant) measure $\mu (A) = \int_A dt/t$ 
for  any measurable set $A \subset \mathbb{R}^+.$ Finally, by $X^\infty_c$ we will denote the space of all functions $f:\mathbb{R}^+ \rightarrow \mathbb{C}$ 
such that $\|(\cdot) f(\cdot)\|_\infty = \sup_{x>0}|x^cf(x)| < +\infty.$

 The Mellin transform of a function $f\in X^1_c$ is defined by (see e.g. \cite{MA},  \cite{BJ1})
$$ M_c[f](s) \equiv [f]^{\wedge}_{M_c} (s) = \int_0^\infty u^{s-1} f(u) du~~~(s=c+ it, t\in \mathbb{R}).$$
Basic properties of the Mellin transform are the following:
$$M_c[af + bg](s) = a M_c[f](s) + bM_c[g](s)~~~(f,g \in X^1_c,~a,b \in \mathbb{R}),$$
$$|M_c[f](s)| \leq \|f\|_{X^1_c}~~(s = c+it).$$
The inverse Mellin transform $M^{-1}_c[g]$ of the function $g \in L^1(\{c\} \times i\mathbb{R})$ is defined by
$$M^{-1}_c[g](x)  := \frac{x^{-c}}{2 \pi}\int_{-\infty}^{+\infty} g(c+it) x^{-it}dt ~~~(x \in \mathbb{R^+}),$$
where by $L^p(\{c\} \times i\mathbb{R})$ for $p \geq 1,$ will mean the space of all functions $g:\{c\} \times i\mathbb{R} \rightarrow \mathbb{C}$ with 
$g(c +i\cdot) \in L^p(\mathbb{R}).$

For $p=2$ the Mellin transform $M_c^2$ of $f \in X^2_c$ is given by (see \cite{BJ2})
$$M_c^2[f](s) \equiv [f]^{\wedge}_{M_c^2} (s) = \limm_{\rho \rightarrow +\infty}~\int_{1/\rho}^\rho f(u) u^{s-1}du~~~(s=c+it),$$
in the sense that
$$\lim_{\rho \rightarrow \infty}\bigg\|M_c^2[f](c+it) - \int_{1/\rho}^\rho f(u) u^{s-1}du\bigg\|_{L^{2}(\{c\}\times i\mathbb{R})} = 0.$$
In this instance, the Mellin transform is norm-preserving in the sense that (see \cite{BJ2})
$$\|g\|_{X^2_c} = \frac{1}{\sqrt{2\pi}}\|[g]^\wedge_{M_c}\|_{L^2(\{c\}\times i\mathbb{R})}.$$
More generally, using the Riesz-Thorin convexity theorem, one may introduce a definition of Mellin transform in $X^p_c$ with $p\in ]1,2[$ in an analogous 
way, i.e., $$\lim_{\rho \rightarrow \infty}\bigg\|M_c^p[f](c+it) - \int_{1/\rho}^\rho f(u) u^{s-1}du\bigg\|_{L^{p'}(\{c\}\times i\mathbb{R})} = 0,$$
where  $p'$ denotes the conjugate exponent of $p.$

Analogously, the inverse Mellin transform $M_c^{-1, 2}$ of a function $g \in L^{2}(\{c\} \times i\mathbb{R})$ is defined as
$$M_c^{-1,2}[f](s) = \limm_{\rho \rightarrow +\infty}~\int_{-\rho}^\rho g(c + iv) v^{-c-iv}dv,$$
in the sense that
$$\lim_{\rho \rightarrow \infty}\bigg\|M_c^{-1,2}[g](c+iv) - \frac{1}{2\pi}\int_{-\rho}^\rho g(c + iv) v^{-c-iv}dv\bigg\|_{X^2_c} = 0.$$
In a similar way one can define the inverse Mellin transform $M_c^{-1,p}$ with $p \in {]}1,2{[}.$

In what follows, we will continue to denote the Mellin transform of a function $g\in L^p(\mathbb{R}^+)$ by $[g]^\wedge_{M_c}$ and we will consider essentially the cases $p=1$ and $p=2.$ 

The Mellin translation operator $\tau_h^c$ for $h \in \mathbb{R}^+,~c \in \mathbb{R}$ and $f: \mathbb{R}^+ \rightarrow \mathbb{C}$ is defined by
$$(\tau_h^c f)(x) := h^c f(hx)~~(x\in \mathbb{R}^+).$$
\noindent
Setting $\tau_h:= \tau^0_h,$ then  $(\tau_h^cf)(x) = h^c (\tau_hf)(x)$  and
$\|\tau_h^c f\|_{X^1_c} = \|f\|_{X^1_c}.$

\vskip0,4cm
For $1\leq p \leq 2$, denote by $B^p_{c,\sigma}$ the Bernstein space of all functions $f\in X^p_c\cap C(\mathbb{R}^+),$ $c \in \mathbb{R},$ which are Mellin 
band-limited 
to the interval $[-\sigma,\sigma],$ $\sigma \in \mathbb{R}^+,$ thus for which $[f]^\wedge_{M_c}(c+it) = 0$ for all $|t| > \sigma.$ We notice that, as in Fourier analysis, the inclusion $B^{p}_{c,\sigma} \subset B^q_{c,\sigma}$ holds for $1 \leq p<q \leq 2.$


\section{A notion of distance}

For $q \in [1,+\infty]$, let $G_c^q$ be the linear space of all functions $f:\mathbb{R}^+ \rightarrow \mathbb{C}$ that have the representation
$$f(x) = \frac{1}{2\pi}\int_{-\infty}^{\infty}\varphi(v) x^{-c-iv}dv \qquad (x>0),$$
where $\varphi \in L^1(\mathbb{R}) \cap L^q(\mathbb{R}).$

The space $G_c^q$ will be endowed with the norm
$$[\!\!|f|\!\!]_q := \|\varphi\|_{L^q(\mathbb{R})} = \left(\int_{-\infty}^\infty |\varphi(v)|^qdv\right)^{1/q}.$$
Note that this is really a norm. Indeed, $[\!\!|f|\!\!]_q = 0$ iff $f=0$ due to the existence and uniqueness of Mellin inversion (see \cite{BJ1}, \cite{BJ2}).

The above norm induces the metric
$$\mbox{\rm dist}_q(f,g) := [\!\!|f-g|\!\!]_q \quad \quad f,g \in G_c^q.$$
Note that in case $q=2$ we have
$$\mbox{\rm dist}_2(f,g) = \sqrt{2\pi}\|f-g\|_{X^2_c},$$
i.e., our distance reduces to the ``Euclidean'' distance in $X^2_c,$ up to the factor $\sqrt{2\pi}.$\newline
\noindent
As a consequence of the Mellin inversion formula, functions $f$ for which $[f]^\wedge_{M_c}(c+i\cdot) \in L^1(\mathbb{R}) \cap L^q(\mathbb{R})$ belong 
to $G_c^q.$  For $p \in [1,2]$, the Mellin-Bernstein space $B_{c,\sigma}^p$ is a subspace of $G_c^q$ since the Mellin transform of $f \in B_{c,\sigma}^p$ 
has compact support as a function of $v \in \mathbb{R}$ and so it belongs to any $L^q(\mathbb{R}).$

For $f \in G_c^q$ we define
$$\mbox{\rm dist}_q(f,B_{c,\sigma}^p) = \inf_{g \in B_{c,\sigma}^p}[\!\!|f-g|\!\!]_q\,.$$
\vskip0,4cm

 The following representation theorem holds:
\newtheorem{Theorem}{Theorem}
\begin{Theorem} \label{representation1} For any $f \in G_c^q$, we have
$$\mbox{\rm dist}_q(f,B_{c,\sigma}^p) = \left(\int_{|v| \geq \sigma}|\varphi (v)|^qdv\right)^{1/q} \quad \quad (1\leq q<\infty),$$
and if $\varphi$ is also continuous, then
$$\mbox{\rm dist}_\infty(f,B_{c,\sigma}^p) = \sup_{|v| \geq \sigma}|\varphi(v)|.$$
\end{Theorem}
{\bf Proof}. 
Assume $q< \infty.$ Clearly, since $g \in B_{c, \sigma}^p$ implies $|[g]^\wedge_{M_c}(c+iv)| = 0$ for $|v| \geq \sigma,$ one has
$$\mbox{\rm dist}_q(f,B_{c,\sigma}^p) = \left\{\inf_{g \in B_{c,\sigma}^p}\int_{|v|\leq \sigma}|\varphi(v) - [g]^\wedge_{M_c}(c+iv)|^qdv + \int_{|v|\geq \sigma}|\varphi(v)|^qdv\right\}^{1/q}.$$
Therefore we have to prove that
$$I_{p,q}:= \inf_{g \in B_{c,\sigma}^p}\int_{|v| \leq \sigma} |\varphi(v) - [g]^\wedge_{M_c}(c+iv)|^q dv = 0.$$
For the given $\sigma >0,$ the space $C^\infty_c({]}-\sigma, \sigma{[})$, whose elements are all the infinitely differentiable functions with compact 
support in ${]}-\sigma, \sigma{[},$ is dense in $L^q({]}-\sigma, \sigma{[})$ for $1\leq q< \infty$ (see e.g. \cite{AD}). Thus, given $\varphi 
\in L^q(\mathbb{R})$ and $\varepsilon >0$, we can take a function $P \in C^\infty_c({]}-\sigma, \sigma{[})$ such that 
$$\|\varphi - P\|_{L^q({]}-\sigma, \sigma{[})} < \varepsilon.$$
Now we define
$$g_\varepsilon (x) := \frac{x^{-c}}{2\pi}\int_{-\sigma}^\sigma P(v) x^{-iv}dv \quad (x>0).$$
Integrating $k$-times by parts, one can easily see that $x^cg_\varepsilon(x) = {\cal O}((\log x)^{-k})$ for $x \rightarrow +\infty$ and $x\rightarrow 0^+.$ 
This implies that $g_\varepsilon \in X^p_c$ for any $p\geq 1$, and $[g_\varepsilon]^\wedge_{M_c}(c+iv) = P(v)$ for $|v| \leq \sigma$ and $0$ otherwise. 
Thus $g_\varepsilon \in B^p_{c, \sigma}.$ Now we conclude that
\begin{eqnarray*}
I_{p,q}^{1/q} &\leq& \|\varphi - [g_\varepsilon]^\wedge_{M_c}(c+i\cdot)\|_{L^q(]-\sigma, \sigma[)} \\&=& \|\varphi - P\|_{L^q(]-\sigma, \sigma[)} < \varepsilon.
\end{eqnarray*}
Hence the assertion follows for $1\leq q < \infty.$

The case $q=\infty$ is treated in a different way. To this end, given $\varepsilon >0,$ there  exists a twice continuously differentiable function $\psi$ 
on $\mathbb{R}$ such that $\sup_{v \in \mathbb{R}}|\varphi (v) - \psi (v)| < \varepsilon/2.$ For example, $\psi$ may be chosen as an appropriate spline 
function. For $\eta \in {]}0,\sigma{[},$ define
\begin{eqnarray*}
\psi_1(x):= \left\{\begin{array}{llll} 0 \quad &\mbox{if}\quad |v|\leq \sigma -\eta,\\[2ex]
\frac{\psi(-\sigma)}{-\eta^3}(v+\sigma-\eta)^3 \quad &\mbox{if}\quad -\sigma \leq v \leq -\sigma +\eta,\\[2ex] 
\frac{\psi(\sigma)}{\eta^3}(v-\sigma + \eta)^3 \quad &\mbox{if}\quad \sigma - \eta \leq v \leq \sigma,\\[2ex]
\psi(v) \quad &\mbox{if}\quad |v| \geq \sigma. 
\end{array} \right.
\end{eqnarray*}
Note that $\psi_1$ is continuous on $\mathbb{R}$ and 
$$\|\psi_1\|_{L^\infty(\mathbb{R})} = \sup_{|v| \geq \sigma -\eta}|\psi_1(v)| = \sup_{|v| \geq \sigma}|\psi (v)| \leq
\sup_{|v|\geq \sigma}|\varphi(v)| + \frac{\varepsilon}{2}\,.$$
Next define $\psi_0(v):= \psi(v) - \psi_1(v)$ for $v \in \mathbb{R}.$ Then $\psi_0$ is continuous on $\mathbb{R},$ twice continuously differentiable on 
${]}-\sigma, \sigma{[},$ it vanishes at $\pm \sigma$ and it has support on $[-\sigma, \sigma].$ With these properties, two integrations by parts show that
$$g_\varepsilon(x):=\frac{x^{-c}}{2 \pi}\int_{-\sigma}^\sigma \psi_0(v)x^{-iv}dv \quad \quad (x>0)$$
defines a function $g_\varepsilon \in X^p_c \cap B^p_{c,\sigma}.$ Furthermore, the Mellin inversion formula yields that 
$[g_\varepsilon]^\wedge_{M_c}(c+i\cdot) = \psi_0(\cdot).$ Now we conclude that
\begin{eqnarray*}
\mbox{\rm dist}_\infty(f, B^p_{c,\sigma}) &= & \inf_{g \in B^p_{c, \sigma}}\|\varphi - [g]^\wedge_{M_c}(c+i\cdot)\|_{L^\infty(\mathbb{R})}\\&\leq& 
\|\varphi - [g_\varepsilon]^\wedge_{M_c}(c+i\cdot)\|_{L^\infty(\mathbb{R})}\\
&=&\|\varphi - \psi_0\|_{L^\infty(\mathbb{R})}\\
&\leq& \|\varphi - \psi\|_{L^\infty(\mathbb{R})} +\|\psi - \psi_0\|_{L^\infty(\mathbb{R})}\\ &\leq& 
\frac{\varepsilon}{2} + \|\psi_1\|_{L^\infty(\mathbb{R})}\\
&\leq& \varepsilon + \sup_{|v| \geq \sigma} |\varphi (v)|.
\end{eqnarray*}
This implies that $\mbox{\rm dist}_\infty(f, B^p_{c, \sigma}) \leq \sup_{|v| \geq \sigma}|\varphi (v)|.$ On the other hand,
\begin{eqnarray*}
\mbox{\rm dist}_\infty(f, B^p_{c, \sigma}) &=& \max\left\{\inf_{g \in B^p_{c, \sigma}} \sup_{|v| \leq \sigma}|\varphi (v) - [g]^\wedge_{M_c}(c+iv)|, \sup_{|v| \geq \sigma}|\varphi (v)|\right\}\\&\geq& \sup_{|v| \geq \sigma}|\varphi(v)|.
\end{eqnarray*}
Hence the formula stated in the theorem holds. \hfill$\Box$

\vskip0,4cm
Next we will obtain a distance formula for Mellin derivatives.
We define the first Mellin derivative (the Mellin differential operator of  first order) by
$$(\Theta_c^1f)(x) := \lim_{h\rightarrow 1}\frac{\tau^c_hf(x) - f(x)}{h-1} =
\lim_{h\rightarrow 1}\frac{h^cf(hx) - f(x)}{h-1} \quad (x>0),$$
and the Mellin differential operator of order $r \in \mathbb{N}$ is defined iteratively by $\Theta^r_c:= \Theta_c^1(\Theta_c^{(r-1)});$ see \cite{BJ1}. 
We have the following
\begin{Theorem}\label{derivative} Let $f\in G_c^q.$ If $v^r\varphi (v)$ belongs to $L^1(\mathbb{R})$ as a function of $v$ for some $r \in \mathbb{N},$ 
then $f$ has Mellin derivatives up to order $r$ in $C_0(\mathbb{R})$ and 
$$(\Theta^k_cf)(x) = \frac{(-i)^k}{2 \pi}\int_{-\infty}^{+\infty} v^k\varphi (v)x^{-c-iv}dv  \quad \quad(k=0,1,\ldots r).$$
\end{Theorem}
{\bf Proof.} Suppose that $r=1.$ For $h \neq 1$ we have
\begin{eqnarray*}
\frac{h^cf(hx) - f(x)}{h-1} &=& \frac{1}{2 \pi}\frac{1}{h-1}\bigg(\int_{-\infty}^{+\infty} \varphi(v) h^{-iv}x^{-c-iv}dv - 
\int_{-\infty}^{+\infty} \varphi(v) x^{-c-iv}dv\bigg) \\[1ex]
&=& \frac{x^{-c}}{2 \pi}\int_{-\infty}^{+\infty} \varphi(v) x^{-iv}\frac{h^{-iv} - 1}{h - 1}dv.
\end{eqnarray*}
Now
\begin{eqnarray*}
\bigg|\frac{h^{-iv} - 1}{h - 1}\bigg|= 
\frac{2}{|h-1|}\bigg|\frac{e^{-\frac{iv \log h}{2}}-e^{\frac{iv \log h}{2}}}{2i}\bigg|=
2 \bigg|\frac{\sin (\frac{v \log h}{2})}{h-1}\bigg| \leq |v|.
\end{eqnarray*}
Since 
$$\lim_{h \rightarrow 1} \frac{h^{-iv} - 1}{h - 1} = -iv$$
and $v\varphi(v)$ is absolutely integrable, Lebesgue's theorem on dominated convergence gives
$$(\Theta^1_c f)(x) = \frac{-i}{2\pi}\int_{-\infty}^{+\infty} v\varphi(v) x^{-c-iv}dv.$$
Moreover, using the Mellin inversion theorem (see \cite[Lemma~4, p.~349]{BJ1},
we have $\Theta^1_cf \in C_0(\mathbb{R}^+).$

The proof for general $r$ follows by mathematical induction. \hfill$\Box$
\vskip0,3cm
Using Theorems \ref{representation1} and \ref{derivative}, we obtain immediately the distance of $\Theta^k_cf$ from the Bernstein space $B_{c,\sigma}^p$.
\vskip0,4cm
\noindent
\newtheorem{Corollary}{Corollary}
\begin{Corollary}\label{cor1}
Let $f \in G^q_c$ with $v^k\varphi \in L^1(\mathbb{R}) \cap L^q(\mathbb{R}).$ Then for every $p \in [1,2]$, we have
$$\mbox{\rm dist}_q(\Theta^k_cf,B_{c,\sigma}^p) = \left(\int_{|v| \geq \sigma}|v^k\varphi (v)|^qdv\right)^{1/q} \quad \quad (1\leq q<\infty),$$
and if $\varphi$ is continuous, then
$$\mbox{\rm dist}_\infty(\Theta^k_cf,B_{c,\sigma}^p) = \sup_{|v| \geq \sigma}|v^k\varphi(v)|.$$
\end{Corollary}
\vskip0,4cm


\section{Estimation of the Mellin distance}

In this section we will introduce certain basic ``intermediate'' function spaces between the spaces $B^p_{c,\sigma}$ and the space $G_c^q.$ 
We will consider mainly the cases $p=1$ and $p=2.$ 
In the following for $p\in [1,2],$ we will denote by $\mathcal{M}_c^p$ the space comprising all functions $f \in X_c^p \cap C(\mathbb{R})$ such that 
$[f]^\wedge_{M_c} \in L^1(\{c\}\times i \mathbb{R}).$ This space is contained in $G_c^q$ for suitable values of $q,$ namely for
$q \in [1, p']$ with $p'$ being the conjugate exponent of $p.$ As for the classes $B^p_{c,\sigma},$ we have again the inclusion $\mathcal{M}^p_c \subset \mathcal{M}^q_c$ for $1 \leq p<q\leq 2.$

We begin with the definitions of differences of integer order and an appropriate modulus of smoothness.
For a function $f \in X^p_c,$ $r \in \mathbb{N}$ and $h >0$, we define
$$(\Delta_h^{r,c}f)(u) := \sum_{j=0}^r (-1)^{r-j}\left(\begin{array}{ll} r \\ j \end{array}\right) f(h^ju)h^{jc},$$
and for $\delta >0,$
$$\omega_{r}(f, \delta, X^p_c) := \sup_{|\log h|\leq \delta}\|\Delta_h^{r,c}f\|_{X^p_c}.$$
In particular for $p=1,$ 
$$\omega_{r}(f, \delta, X^1_c) := \sup_{|\log h|\leq \delta}\|\Delta_h^{r,c}f\|_{X^1_c}.$$
Among the basic properties of the above modulus of smoothness $\omega_{r}$ we list the following three:
\begin{enumerate}
\item $\omega_r(f, \cdot, X^p_c)$ is a non decreasing function on $\mathbb{R}^+;$
\item  $\omega_r(f, \delta, X^p_c) \leq 2^r\|f\|_{X^p_c}$;
\item for any  positive $\lambda$ and $\delta$, one has
$$\omega_r(f, \lambda \delta, X^p_c) \leq (1+\lambda)^r\omega_r(f,\delta,X^p_c).$$
\end{enumerate}

We know that for functions $f \in X^1_c$ or $f \in X^2_c$ one has (see \cite{BJ1},  \cite{BJ0}, \cite{BBM})
\begin{eqnarray}
[\Delta_h^{r,c}f]^\wedge_{M_c}(c+iv) = (h^{-iv} - 1)^r [f]^\wedge_{M_c}(c+iv) \qquad (v \in \mathbb{R}).
\end{eqnarray}
We have the following
\begin{Theorem} \label{estimate1} 
 If $f \in \mathcal {M}^1_c,$ then for any $q \in [1, \infty]$, 
\begin{eqnarray*}
\mbox{\rm dist}_q(f, B^1_{c, \sigma}) \leq D \cdot \left\{\begin{array}{ll}
\displaystyle\left\{\int_{\sigma}^\infty [\omega_r(f, v^{-1}, X^1_c)]^qdv\right\}^{1/q} &\quad (q< \infty),\\[2ex] 
\omega_r(f, \sigma^{-1}, X^1_c) &\quad (q=\infty),
\end{array} \right.
\end{eqnarray*}
where $D$ is a constant depending on $r$ and $q$ only.
\end{Theorem}
{\bf Proof}. From (1), setting $h= e^{\pi/v}$, we have
$$[\Delta_{h}^{r,c}f]^\wedge_{M_c}(c+iv) = (-2)^r [f]^\wedge_{M_c}(c+iv)$$
or
$$[f]^\wedge_{M_c}(c+iv) = \frac{1}{(-2)^r}\int_0^\infty (\Delta^{r,c}_{h}f)(u) u^{c+iv-1}du \quad \quad (h= e^{\pi/v})$$
and so
$$|[f]^\wedge_{M_c}(c+iv)| \leq \frac{1}{2^r}\int_0^\infty |(\Delta^{r,c}_{h}f)(u)|u^{c-1}du \leq \frac{1}{2^r}\omega_r(f, \frac{\pi}{|v|}, X^1_c).$$
Now using the properties of the modulus $\omega_r,$ we find that
$$\omega_r(f, \frac{\pi}{|v|}, X^1_c) \leq (1+\pi)^r \omega_r(f, \frac{1}{|v|}, X^1_c).$$
Thus
$$|[f]^\wedge_{M_c}(c+iv)| \leq \left(\frac{1+\pi}{2}\right)^r \omega_r(f, \frac{1}{|v|}, X^1_c).$$
In view of Theorem~\ref{representation1}, this implies the assertion for $q<\infty$.
The case $q=\infty$ is obtained analogously. \hfill $\Box$
\vskip0,3cm
\begin{Theorem} \label{estimate2} 
If $f \in \mathcal{M}^2_c$,  then for any $q\in [1,2],$
$$\mbox{\rm dist}_q(f, B^2_{c,\sigma}) \leq D \left\{\int_{\sigma}^\infty [v^{-q/2}\omega_r(f, v^{-1}, X^2_c)]^qdv\right\}^{1/q},$$
where $D$ is a constant depending on $r$ and $q$ only.
\end{Theorem}
{\bf Proof}. First we consider the case $q=2.$ Then 
$$|[\Delta_h^{r,c}f]^\wedge_{M_c}(c+iv)| = 2^r |\sin ((v \log h)/2)| |[f]^\wedge_{M_c}(c+iv)|$$
and (see \cite[Lemma 2.6]{BJ2})
$$\|[\Delta_h^{r,c}f]^\wedge_{M_c}\|_{L^2(\{c\}\times i \mathbb{R})} = \sqrt{2\pi}\|\Delta_h^{r,c}f\|_{X^2_c} \leq \sqrt{2\pi}\omega_r(f, |\log h|, X^2_c).$$
Now let $h\geq 1.$ For $v \in [(2\log h)^{-1}, (\log h)^{-1}]$, one has 
$$\sin ((v\log h)/2) \geq \frac{1}{2\pi}\,,$$
and hence 
$$\int_{1/(2\log h)}^{1/\log h} |[f]^\wedge_{M_c}(c+iv)|^2 dv \leq 
(2\pi)^{2r}\int_0^\infty |\sin ((v\log h)/2)|^{2r}| |[f]^\wedge_{M_c}(c+iv)|^2dv.$$
Analogously we have
$$\int_{-1/\log h}^{-1/(2\log h)} |[f]^\wedge_{M_c}(c+iv)|^2 dv \leq 
(2\pi)^{2r}\int_{-\infty}^0 |\sin ((v\log h)/2)|^{2r}| |[f]^\wedge_{M_c}(c+iv)|^2dv,$$
and so
\begin{eqnarray*}
\lefteqn{\int_{1/(2\log h) \leq|v| \leq 1/\log h}|[f]^\wedge_{M_c}(c+iv)|^2 dv} \qquad\qquad\quad\\
&\leq &
(2\pi)^{2r}\int_{-\infty}^\infty |\sin ((v\log h)/2)|^{2r}| |[f]^\wedge_{M_c}(c+iv)|^2dv \\
&\leq& \pi^{2r}[\omega_r(f, \log h, X^2_c)]^2.
\end{eqnarray*}
Now, let $\sigma >0$ be fixed, set $\sigma_k:= \sigma 2^k$ with $k \in \mathbb{N}_0$, and define $h$ by $\log h = 1/\sigma_{k+1}.$ Then
$$\int_{\sigma_k \leq |v| \leq \sigma_{k+1}}|[f]^\wedge_{M_c}(c+iv)|^2dv \leq 
\pi^{2r}[\omega_r(f, \sigma_{k+1}^{-1}, X^2_c)]^2,$$
and so summation over $k$ yields
$$\int_{|v| \geq \sigma}|[f]^\wedge_{M_c}(c+iv)|^2dv \leq 
\pi^{2r}	\sum_{k=0}^\infty [\omega_r(f, \sigma_{k+1}^{-1}, X^2_c)]^2.$$
Now since $\sigma_{k+1}- \sigma_k = \sigma_k,$ from the monotonicity of $\omega_r$ as a function of $\delta,$ one has
$$\int_{\sigma_k}^{\sigma_{k+1}}v^{-1}[\omega_r(f, v^{-1}, X^2_c)]^2dv \geq 
\frac{\sigma_k}{\sigma_{k+1}}[\omega_r(f, \sigma^{-1}_{k+1}, X^2_c)]^2,$$
from which we deduce
$$\sum_{k=0}^\infty [\omega_r(f, \sigma^{-1}_{k+1}, X^2_c)]^2 \leq 
2\int_\sigma^\infty v^{-1}[\omega_r(f, v^{-1}, X^2_c)]^2dv.$$
This gives the assertion  for $q=2.$ For $q\in [1, 2{[}$ one can proceed as in the proof of Proposition 13 in \cite{BSS2}, using H\"{o}lder's 
inequality. \hfill $\Box$


\subsection{Mellin-Lipschitz spaces}

For $\alpha \in {]}0,r]$ we define the Lipschitz class by
$$\mbox{Lip}_r(\alpha, X^p_c) := \{f \in X^p_c: \omega_r(f;\delta;X^p_c) = {\cal O}(\delta^\alpha), \delta \rightarrow 0^+\}.$$
\noindent
As a consequence of Theorems \ref{estimate1} and \ref{estimate2}, we obtain the following corollary which determines the Mellin distance of a function 
$f \in \mbox{Lip}_r(\alpha, X^p_c)$ from $B^p_{c,\sigma}$ for  $p=1,2.$
\begin{Corollary}\label{lip}
If $f \in \mbox{\rm Lip}_r(\alpha, X^1_c\cap C(\mathbb{R}))$ for some $r\in \mathbb{N},~ r \geq 2$ and $1<\alpha \leq r,$ then
$$\mbox{\rm dist}_1(f, B^{1}_{c, \sigma}) = {\cal O}(\sigma^{-\alpha +1})\quad \quad (\sigma \rightarrow +\infty).$$
Moreover, if $f \in \mathcal{M}^2_c \cap \mbox{\rm Lip}_r(\beta, X^2_c)$ with $r \in \mathbb{N},$ $q^{-1} - 2^{-1} < \beta \leq r,$ then
$$\mbox{\rm dist}_q (f, B^2_{c, \sigma}) = {\cal O}(\sigma^{-\beta - 1/2 + 1/q})\quad \quad (\sigma \rightarrow +\infty).$$
\end{Corollary}
The proof follows immediately from Theorems \ref{estimate1} and \ref{estimate2}.  Note that, if  $f \in \mbox{Lip}_r(\alpha, X^1_c\cap C(\mathbb{R})),$ 
then from the proof of Theorem \ref{estimate1} one has that $[f]^\wedge_{M_c} \in L^1(\{c\}\times i\mathbb{R});$ 
thus $f \in\mathcal{M}^1_c.$ 

For $q=2$ we obtain the estimate
$$\mbox{\rm dist}_2 (f, B^2_{c, \sigma}) = {\cal O}(\sigma^{-\beta})\quad \quad (\sigma \rightarrow +\infty).$$


\subsection{Mellin-Sobolev spaces}

Denote by $AC_{{\tt loc}}(\mathbb{R}^+)$ the space of all locally absolutely continuous functions on $\mathbb{R}^+$. The Mellin-Sobolev space 
$W_c^{r,p}(\mathbb{R}^+)$ is defined as the space of all functions $f \in X^p_c$ which are equivalent to a function 
$g \in C^{r-1}(\mathbb{R}^+)$ with 
$g^{(r-1)}\in AC_{{\tt loc}}(\mathbb{R}^+)$ such that $\Theta^r_cg \in X^p_c$ (see \cite{BJ1}, \cite{BJ2}, \cite{BBM}). 
For $p=1$ it is well known that for any $f \in W_c^{r,1}$ one has (see \cite{BJ1})
$$[\Theta^r_c]^\wedge_{M_c}(c+iv) = (-iv)^r[f]^\wedge_{M_c}(c+iv)\quad \quad (v \in \mathbb{R}).$$
The same result also holds for $1< p\leq 2,$ taking into account the general convolution theorem for Mellin transforms 
(\cite[Lemma~3.1]{BJ2} in case $p=2$, \cite[Lemma 2]{BKT}). 

By the above result, the Mellin-Sobolev space $W_c^{r,p}(\mathbb{R}^+)$ can be characterized as
$$W_c^{r,p}(\mathbb{R}^+) = \{f \in X^p_c: (-iv)^r [f]^\wedge_{M_c}(c+iv) = [g]^\wedge_{M_c}(c+iv), ~g \in L^p(\{c\}\times i \mathbb{R})\}.$$
We have the following
\begin{Theorem}\label{sobolev1}
Let $f \in \mathcal{M}^1_c\cap W^{r,1}_c(\mathbb{R}^+).$ Then for $q \in [1,\infty]$ and $r >1/q,$
\begin{eqnarray*}
\mbox{\rm dist}_q(f, B^1_{c,\sigma}) \leq D\|\Theta^r_cf\|_{X^1_c} \cdot \left\{\begin{array}{ll} \sigma^{-r +1/q}, &\quad q<\infty,\\ 
\sigma^{-r}, &\quad q= \infty,
\end{array} \right.
\end{eqnarray*}
where $D$ is a constant depending on $r$ and $q$ only.
If, in addition, $v[f]^\wedge_{M_c}(c+iv) \in L^1(\mathbb{R})$, then for $r> 1 + 1/q,$
\begin{eqnarray*}
\mbox{\rm dist}_q(\Theta_c f, B^1_{c,\sigma}) \leq D'\|\Theta^r_cf\|_{X^1_c} \cdot \left\{\begin{array}{ll} \sigma^{-r +1+1/q},  & \quad q<\infty,\\ 
\sigma^{-r+1}, & \quad q= \infty,
\end{array} \right.
\end{eqnarray*}
where $D'$ is again a constant depending on $r$ and $q$ only.
\end{Theorem}
{\bf Proof}. First we consider the case $q<\infty.$ 
 The formula for the Mellin transform of a Mellin derivative yields
$$[f]^\wedge_{M_c}(c+iv) = (-iv)^{-r}[\Theta^r_cf]^\wedge_{M_c}(c+iv) \qquad (v \in \mathbb{R}
\setminus \{0\}).$$
Thus, from Theorem \ref{representation1} we obtain 
\begin{eqnarray*}
\mbox{\rm dist}_q(f, B^1_{c,\sigma}) = \left\{\int_{|v| \geq \sigma}|v^{-r}[\Theta^r_cf(v)]^\wedge_{M_c}(c+iv)|^qdv\right\}^{1/q}.
\end{eqnarray*}
    Since $\Theta^r_cf \in X^1_c,$ its Mellin transform is continuous and bounded on $\{c\} \times i\mathbb{R}$ (see \cite{BJ1}). Therefore 
\begin{eqnarray*}
\mbox{\rm dist}_q(f, B^1_{c,\sigma}) &\leq& \|[\Theta^r_cf]^\wedge_{M_c}\|_{C(\{c\}\times i \mathbb{R})} \left\{2\int_{v \geq 
\sigma}v^{-rq}dv\right\}^{1/q} \\
&\leq& \|\Theta^r_cf\|_{X^1_c}\bigg(\frac{2}{rq-1}\bigg)^{1/q} \frac{1}{\sigma^{r -1/q}},
\end{eqnarray*}
and hence the assertion for $q<\infty$ is proved with $D= (2/(rq-1))^{1/q}.$
For $q= \infty$ we use again Theorem \ref{representation1} and proceed analogously, obtaining $D=1.$ For the second part, note that under the assumptions 
on $f$ and $v [f]^\wedge_{M_c}(c+iv) \in L^1(\{c\}\times i \mathbb{R})$, we have $\Theta_cf \in {M}^1_c \cap W^{r-1, 1}_c(\mathbb{R}^+).$ 
Therefore we can apply the first part of the proof to the function $\Theta_cf,$ obtaining immediately the assertion with the constant 
$D'= (2/(rq-q-1))^{1/q}$ for $q< \infty$ and $D'= 1$ for $q = \infty.$
\hfill$\Box$
\vskip0,3cm 
Note that if $f\in \mathcal{M}^1_c\cap W^{r,1}(\mathbb{R}^+)$ satisfies the further condition that 
$[\Theta^r_cf]^\wedge_{M_c} \in L^q(\{c\} \times i\mathbb{R}),$ then one may write
\begin{eqnarray*}
\mbox{\rm dist}_q(f, B^1_{c,\sigma}) &=& \left\{\int_{|v| \geq \sigma}|v^{-r}[\Theta^r_cf(v)]^\wedge_{M_c}(c+iv)|^qdv\right\}^{1/q}\\[2ex]
&\leq& \frac{1}{\sigma^r}\|[\Theta^r_cf]^\wedge_{M_c}\|_{L^q(\{c\} \times i\mathbb{R})}.
\end{eqnarray*}
Moreover, one has 
$$\mbox{\rm dist}_q(f, B^1_{c,\sigma}) = \mathcal{O}(\sigma^{-r}) \qquad (\sigma \rightarrow +\infty).$$
\vskip0,4cm

For $p=2$ we have the following
\begin{Theorem}\label{sobolev2}
Let $f \in \mathcal{M}^2_c\cap W^{r,2}_c(\mathbb{R}^+).$ Then for $q \in [1,2]$, 
\begin{eqnarray*}
\mbox{\rm dist}_q(f, B^2_{c,\sigma}) \leq D\|\Theta^r_cf\|_{X^2_c}~ \sigma^{-r -1/2 + 1/q},
\end{eqnarray*}
where $D$ is a constant depending on $r$ and $q$ only.
If, in addition, $v[f]^\wedge_{M_c}(c+iv) \in L^1(\mathbb{R})$, then for $r> 1 + 1/2 +1/q,$
\begin{eqnarray*}
\mbox{\rm dist}_q(\Theta_c f, B^2_{c,\sigma}) \leq D'\|\Theta^r_cf\|_{X^2_c}~ \sigma^{-r +1/2+1/q},  
\end{eqnarray*}
where $D'$ is again a  constant depending on $r$ and $q$ only.
\end{Theorem}
{\bf Proof}. As in the previous theorem, by the formula of Mellin transform in $X^2_c$ for derivatives we have
$$\mbox{\rm dist}_q(f, B^2_{c,\sigma}) = \left\{\int_{|v| \geq \sigma}|v^{-r}[\Theta^r_c]^\wedge_{M_c}(c+iv)|^qdv\right\}^{1/q}.$$
For $q=2,$ using the property that the Mellin transform in $X^2_c$ is norm-preserving (see \cite[Lemma~2.6]{BJ2}), we have
\begin{eqnarray*}
\mbox{\rm dist}_q(f, B^2_{c,\sigma})&\leq& \frac{1}{\sigma^r}\left\{\int_{|v| \geq \sigma}|[\Theta^r_c]^\wedge_{M_c}(c+iv)|^2dv\right\}^{1/2} \\[1ex]
&\leq& \frac{1}{\sigma^r}\|[\Theta^r_cf]^\wedge_{M_c}\|_{L^2(\{c\}\times i \mathbb{R})} = \sqrt{2\pi}\frac{1}{\sigma^r}\|\Theta^r_cf\|_{X^2_c}.
\end{eqnarray*}
Therefore the assertion follows for $q=2$ with the constant $D = (2\pi)^{-1/2}.$ For $q\in [1,2{[}$ one can use  H\"{o}lder's inequality with 
$\mu = 2/(2-q),~\nu = 2/q$, obtaining
\begin{eqnarray*}
\mbox{\rm dist}_q(f, B^2_{c,\sigma})&\leq& 
\left\{2 \int_{\sigma}^\infty v^{-rq\mu}dv\right\}^{1/(q\mu)}
\left\{\int_{|v|\geq \sigma}|[\Theta^r_cf]^\wedge_{M_c}(c+iv)|^{q\nu}dv\right\}^{1/(q\nu)}\\[1ex] &\leq&
\frac{1}{\sigma^{r+1/2 -1/q}}\left\{\frac{4-2q}{(2r+1)q - 2}\right\}^{1/q - 1/2} 
\|[\Theta^r_cf]^\wedge_{M_c}\|_{L^2(\{c\}\times i \mathbb{R})} \\[1ex]
&=& \sqrt{2\pi}
\frac{1}{\sigma^{r+1/2 -1/q}}\left\{\frac{4-2q}{(2r+1)q - 2}\right\}^{1/q - 1/2} \|\Theta^r_cf\|_{X^2_c}.
\end{eqnarray*}
Thus the first inequality holds with  
$$D=  \sqrt{2\pi}\left\{\frac{4-2q}{(2r+1)q - 2}\right\}^{1/q - 1/2}.$$
The second inequality follows by arguments similar to those in the proof of Theorem~\ref{sobolev1}. \hfill$\Box$
\vskip0,4cm
As a consequence, under the assumptions of the first part of Theorem \ref{sobolev2}, we can obtain an asymptotic estimate of the form
$$\mbox{\rm dist}_q(f, B^2_{c,\sigma}) = \mathcal{O}(\sigma^{-r-1/2+1/q}) \qquad (\sigma \rightarrow +\infty).$$


\section{Applications}

In this section we will illustrate applications to various basic formulae such as the approximate exponential sampling theorem, the approximate 
reproducing kernel formula in the Mellin frame (see \cite{BJ3}, \cite{BBM0}), a generalized Boas differentiation formula and an extension of a 
Bernstein-type inequality.

In the following for $c \in \mathbb{R}$, we denote by $\mbox{lin}_c$  the function
$$\mbox{lin}_c(x) := \frac{x^{-c}}{2\pi i}\frac{x^{\pi i} -x^{-\pi i}}{\log x} = \frac{x^{-c}}{2\pi}\int_{-\pi}^\pi x^{-it}dt \qquad (x>0, \, x\neq 1)$$
with the continuous extension $\mbox{lin}_c(1) = 1.$ Thus 
$$\mbox{lin}_c(x) = x^{-c}\mbox{sinc}(\log x) \qquad (x>0).$$
Here, as usual, the ``sinc'' function is defined as 
$$\mbox{sinc}(t):= \frac{\sin (\pi t)}{\pi t}~ \mbox{for}~t\neq 0,\qquad \mbox{sinc}(0) = 1.$$
It is clear that $\mbox{lin}_c \not \in X_{\overline{c}}$ for any $\overline{c}.$ However, it belongs to the space $X^2_c$ and its Mellin transform in 
$X^2_c$-sense is given by
$$[\mbox{\rm lin}_c]^\wedge_{M_c}(c+iv) = \chi_{[-\pi, \pi]},$$
where $\chi_A$ denotes the characteristic function of the set $A.$

\subsection{Approximate exponential sampling formula}

For a function  $f \in B^2_{c,\pi T}$ the following exponential sampling formula holds (see \cite{BJ3}, \cite{BJ2}):
$$f(x) = \sum_{k \in \mathbb{Z}}f(e^{k/T})\mbox{\rm lin}_{c/T}(e^{-k}x^T) \qquad (x>0).$$
As an approximate version in the space $\mathcal{M}^2_c$ we have (see  \cite[Theorem~5.5]{BJ2}):
\newtheorem{Proposition}{Proposition}
\begin{Proposition}\label{approsamp}
Let $f \in \mathcal{M}_c^2.$
Then there holds the error estimate
\begin{eqnarray*}
\lefteqn{\bigg|f(x) - \sum_{k=-\infty}^{\infty} f(e^{k/T})\mbox{\rm lin}_{c/T}(e^{-k}x^T)\bigg|}\qquad\qquad\quad\\[2ex]
&\leq& \frac{x^{-c}}{\pi}\int_{|t| > \pi T}| [f]^\wedge_{M_c}(c+it)| dt \qquad(x \in \mathbb{R}^+,~T >0).
\end{eqnarray*}
\end{Proposition}
This estimate can now be given a ``metric interpretation''.
By Theorem \ref{representation1}, the right-hand side may be expressed as 
$$\frac{x^{-c}}{\pi}\mbox{\rm dist}_1(f,B^2_{c, \pi T}).$$
Hence, introducing a remainder $(R_{\pi T} f)(x)$ by writing
\begin{equation}\label{approx_samp}
 f(x)\,=\, \sum_{k\in\mathbb{Z}} f\left(e^{k/T}\right) \mbox{lin}_{c/T}\left(e^{-k}x^T\right)
+ \left(R_{\pi T}f\right)(x),
\end{equation}
we have  by Proposition~\ref{approsamp}
$$
|\left(R_{\pi T} f\right)(x)|\,\leq\, \frac{x^{-c}}{\pi}\,\mbox{dist}_1(f, B_{c,\pi
T}^2) \qquad (x>0),$$
or equivalently,
\begin{eqnarray}\label{exp_samp}
\|R_{\pi T} f\|_{X_c^\infty}\,\leq\, \frac{1}{\pi}\,\mbox{dist}_1(f, B_{c,\pi T}^2).
\end{eqnarray}

This relation is a trivial equality when $f\in B_{c, \pi T}^2$. But
equality can also occur when $f\not\in B_{c,\pi T}^2.$
Indeed, consider the function
$$f(x)\,:=\, x^{-c} \mbox{sinc} (2T\log x -1).$$
By a straight forward calculation, we find that
$$ [f]_{M_c}^\wedge(c+iv)\,=\,
\frac{e^{iv/(2T)}}{2T}\,\mbox{rect} \left(\frac{v}{2T}\right),$$
where rect denotes the rectangle function defined by
\begin{eqnarray*}
\mbox{rect}(x)\,:=\left\{
\begin{array}{ccc}
1 & \hbox{ if } & |x|<\pi,\\
\frac{1}{2} & \hbox{ if } &|x|=\pi,\\
0 & \hbox{ if } & |x|>\pi.
\end{array}
\right.
\end{eqnarray*}
Thus, $f\not\in B_{c,\pi T}^2$ and 
$$\mbox{dist}_1(f, B_{c,\pi T}^2)\,=\, \int_{|v|\geq\pi T}
|[f]_{M_c}^\wedge(c+iv)| dv\,=\, \frac{1}{2T} \int_{|v|\geq \pi T}
\mbox{rect} \left(\frac{v}{2T}\right) dv\,=\, \pi.$$
Furthermore,
$$f(e^{k/T})\,=\, e^{-kc/T} \mbox{sinc} (2k-1)\,=0$$
for all $k\in\mathbb{Z}$. Therefore $(R_{\pi T}f)(x)=f(x)$, which shows that
$$\|R_{\pi T}f\|_{X_c^\infty}\,=\, \sup_{x>0} |\mbox{sinc} (2T\log x
-1)|\,=\,1,$$
and so equality occurs in (\ref{exp_samp}).
\vskip0,3cm
Now, employing Theorem \ref{representation1}, Corollary \ref{lip} and the results on Mellin-Sobolev spaces, one has the following theorem.
\begin{Theorem}\label{expestimates}
For the remainder of the approximate exponential sampling formula (\ref{approx_samp}), the following asymptotic estimates hold:
\begin{enumerate}
\item If $f\in \mbox{\rm Lip}_r(\alpha, X^1_c\cap C(\mathbb{R}^+)),$ $r\in \mathbb{N},\,r\geq 2,\,1<\alpha \leq r,$ then
$$\|(R_{\pi T}f)\|_{X^\infty_c} = \mathcal{O}(T^{-\alpha +1})\qquad (T \rightarrow +\infty).$$
\item If $f \in \mathcal{M}^2_c \cap \mbox{\rm Lip}_r(\beta, X^2_c), \,r\in \mathbb{N}, \,1/2 < \beta \leq r,$ then 
$$\|(R_{\pi T}f)\|_{X^\infty_c} = \mathcal{O}(T^{-\beta +1/2})\qquad (T \rightarrow +\infty).$$
\item If $f \in \mathcal{M}^1_c\cap W^{r,1}_c(\mathbb{R}^+), \, r>1,$ then 
$$\|(R_{\pi T}f)\|_{X^\infty_c} = \mathcal{O}(T^{-r+1}) \qquad (T \rightarrow +\infty).$$
\item If $f \in \mathcal{M}^2_c \cap W^{r,2}_c(\mathbb{R}^+), \, r>1/2,$ then 
$$\|(R_{\pi T}f)\|_{X^\infty_c} = \mathcal{O}(T^{-r+1/2}) \qquad (T \rightarrow +\infty).$$
\end{enumerate}
\end{Theorem}


\subsection{Approximate Mellin reproducing kernel formula}

Another interesting formula is the ``Mellin reproducing kernel formula'' for Mellin band-limited functions 
$f \in B^2_{c, \pi T}$. It reads as (see \cite[Theorems 4 and 5]{BBM0})
$$f(x) =  T\int_0^\infty f(y)\mbox{\rm lin}_{c/T}\bigg(\left(\frac{x}{y}\right)^T\bigg)\frac{dy}{y} \qquad (x>0).$$
An approximate version was established in \cite[Theorem 6]{BBM0} for functions in the class $\mathcal{M}^1_c.$ In the same way we can state a version in $\mathcal{M}^2_c,$ as follows
\begin{Proposition}\label{amrkf}
Let $f\in \mathcal{M}_c^2.$
Then for $x \in \mathbb{R}^+,$ and $T>0,$ there holds
\begin{equation}\label{rep_kernel}
f(x) = T\int_0^\infty f(y)\mbox{\rm lin}_{c/T}\bigg(\left(\frac{x}{y}\right)^T\bigg)\frac{dy}{y} + (R^\ast_{\pi T}f)(x),
\end{equation}
where
$$(R^\ast_{\pi T}f)(x) := \frac{x^{-c}}{2 \pi}\int_{|t|\geq \pi T} [f]^\wedge_{M_c}(c+it) x^{-it} dt.$$
Furthermore, we have the error estimate
$$|(R^\ast_{\pi T}f)(x)| \leq \frac{x^{-c}}{2 \pi}\int_{|t|\geq \pi T} |[f]^\wedge_{M_c}(c+it)|dt.$$
\end{Proposition}
{\bf Proof}. The proof is essentially the same as in \cite[Theorem 6]{BBM0}. Here for the sake of completeness we give some details. First, note that the convolution integral in (\ref{rep_kernel}) exists, using a H\"{o}lder-type inequality. Putting $G(x) := \mbox{lin}_{c/T}(x^T),$ its  $X^2_c$~- Mellin transform is given by $[G]^\wedge_{M_c} = T^{-1}\chi_{[-\pi T, \pi T]}.$  Since $[f]^\wedge_{M_c}\in L^1(\{c\}\times i \mathbb{R}),$  using the Mellin inversion and the Fubini Theorem, one can easily obtain, with the same proof, an $X^2_c$~- extension of the Mellin-Parseval formula for convolutions (see \cite[Theorem 9]{BJ1} for functions in $X^1_c$), obtaining
$$T\int_0^\infty f(y)\mbox{\rm lin}_{c/T}\bigg(\left(\frac{x}{y}\right)^T\bigg)\frac{dy}{y} = \frac{x^{-c}}{2\pi}\int_{-\infty}^{+\infty}[f]^\wedge_{M_c}(c+it) [G]^\wedge_{M_c}(c+it)x^{-it}dt.$$
Therefore, by the Mellin inversion formula we have
\begin{eqnarray*}
&&T\int_0^\infty f(y)\mbox{\rm lin}_{c/T}\bigg(\left(\frac{x}{y}\right)^T\bigg)\frac{dy}{y} = \frac{x^{-c}}{2\pi}\int_{-\pi T}^{\pi T}[f]^\wedge_{M_c}(c+it) x^{-it}dt \\
&=&f(x) - \frac{x^{-c}}{2\pi} \int_{|t| \geq \pi T}[f]^\wedge_{M_c}(c+it) x^{-it}dt
\end{eqnarray*}
that is the assertion.
\vskip0,3cm
As before, employing Theorem \ref{representation1}, one can express the error estimate in terms of the distance, i.e., 
$$|(R^\ast_{\pi T}f)(x)|\leq \frac{x^{-c}}{2\pi}\mbox{\rm dist}_1(f,B^2_{c, \pi T})\qquad (x>0),$$
or equivalently,
\begin{equation}\label{ak_est}
\|R^\ast_{\pi T} f\|_{X_c^\infty}\,\leq\, \frac{1}{2\pi}\,\mbox{dist}_1(f, B_{c,\pi T}^2).
\end{equation}
This is again a sharp inequality. Indeed, consider $f(x):=x^{-c} \hbox{sinc}(2T\log x)$. Then $f$ satisfies the hypotheses of
Proposition~\ref{amrkf}. By a calculation we find that
$\hbox{dist}_1(f, B^2_{c,\pi T})=\pi$ and $\|R^\ast_{\pi T} f\|_{X_c^\infty}=1/2.$ Hence equality occurs in (\ref{ak_est}).

Using the estimates of the distance functional in Mellin-Lipschitz and Mellin-Sobolev spaces, we obtain the following results.
\begin{Theorem}\label{rkfest}
For the remainder of the approximate Mellin reproducing kernel formula (\ref{rep_kernel}), the following asymptotic estimates hold:
\begin{enumerate}
\item If $f\in \mbox{\rm Lip}_r(\alpha, X^1_c\cap C(\mathbb{R}^+)),$ $r\in \mathbb{N},\,r\geq 2,\,1<\alpha \leq r,$ then
$$\|(R^\ast_{\pi T}f)\|_{X^\infty_c} = \mathcal{O}(T^{-\alpha +1})\qquad (T \rightarrow +\infty).$$
\item If $f \in \mathcal{M}^2_c \cap \mbox{\rm Lip}_r(\beta, X^2_c), \,r\in \mathbb{N}, \,1/2 < \beta \leq r,$ then 
$$\|(R^\ast_{\pi T}f)\|_{X^\infty_c} = \mathcal{O}(T^{-\beta +1/2})\qquad (T \rightarrow +\infty).$$
\item If $f \in \mathcal{M}^1_c\cap W^{r,1}_c(\mathbb{R}^+), \, r>1,$ then 
$$\|(R^\ast_{\pi T}f)\|_{X^\infty_c} = \mathcal{O}(T^{-r+1}) \qquad (T \rightarrow +\infty).$$
\item If $f \in \mathcal{M}^2_c \cap W^{r,2}_c(\mathbb{R}^+), \, r>1/2,$ then 
$$\|(R^\ast_{\pi T}f)\|_{X^\infty_c} = \mathcal{O}(T^{-r+1/2}) \qquad (T \rightarrow +\infty).$$
\end{enumerate}
\end{Theorem}


\subsection{A sampling formula for Mellin derivatives}

In the context of Fourier analysis the following differentiation formula
has been considered:
\begin{eqnarray}\label{diff1}
f'(x)\,=\, \frac{4T}{\pi} \sum_{k\in\mathbb{Z}} \frac{(-1)^{k+1}}{(2k-1)^2}
f\left(x+ \frac{2k-1}{2T}\right) \qquad (x\in\mathbb{R}).
\end{eqnarray}
It holds for all entire functions of exponential type $\pi T$ which are
bounded on the real line. In particular, it holds for trigonometric
polynomials of degree at most $\lfloor \pi T\rfloor$, where $\lfloor \pi T\rfloor$ denotes the integral part of $\pi T,$ and in this case the
series on the right-hand side can be reduced to a finite sum.

The formula for trigonometric polynomials was discovered by Marcel
Riesz \cite{RIE} in 1914. Its generalization (\ref{diff1}) is due to 
Boas \cite{BOA}. Some authors refer to (\ref{diff1}) as the {\it generalized
Riesz interpolation formula}, others name it after Boas.
  
Formula (\ref{diff1}) has several interesting applications. It provides a
very short proof of Bernstein's inequality in $L^p(\mathbb{R})$ for all $p\in [1,
\infty]$. Modified by introducing a Gaussian multiplier, it leads to a
stable algorithm of high precision for numerical differentiation (see
\cite{SCH2}). Furthermore, it has been extended to higher order
derivatives (see \cite{BSS1}, \cite{SCH2}).

The following theorem gives an analogue of (\ref{diff1}) for Mellin
derivatives.

\begin{Theorem}\label{d_thm1}
For $f\in B_{c,\pi T}^\infty$ there holds
\begin{eqnarray}\label{d_thm1.1}
\Theta_cf(x)\,=\, \frac{4T}{\pi} \sum_{k\in\mathbb{Z}} \frac{(-1)^{k+1}}{(2k-1)^2}
\,e^{(k-1/2)c/T} f\left(x e^{(k-1/2)/T}\right)\qquad (x\in\mathbb{R}^+).
\end{eqnarray}
\end{Theorem}
{\bf Proof}. \,
Formula (\ref{d_thm1.1}) could be deduced from (\ref{diff1}) by making use of
the relationship between the Fourier transform and the Mellin transform.
In the following we give an independent proof completely within Mellin
analysis.

First assume that, in addition,
\begin{eqnarray} \label{d_thm1p1}
x^c f(x)\,=\, \mathcal{O}\left(\frac{1}{|\log x|}\right)
\qquad (x\rightarrow 0_+ \hbox{ and } x\rightarrow \infty).
\end{eqnarray}
Then the exponential sampling formula applies to $f$ and yields
$$ f(x)\,=\, \sum_{k\in\mathbb{Z}} f\left(e^{k/T}\right) \mbox{lin}_{c/T}\left(e^{-k}x^T\right).$$
The series converges absolutely and uniformly on compact subsets of
$\mathbb{R}^+$. When we apply the differentiation operator
$\Theta_c$ with respect to $x$, we may interchange it with the summation on the
right-hand side. Thus,
\begin{eqnarray}\label{d_thm1p2}
\Theta_c f(x)\,=\, \sum_{k\in\mathbb{Z}} f\left(e^{k/T}\right) \Theta_c
\mbox{lin}_{c/T}\left(e^{-k}x^T\right).
\end{eqnarray}
By a calculation we find that
$$\Theta_c\mbox{lin}_{c/T}\left(e^{-k}x^T\right)\,=\, T\,
\frac{e^{kc/T}x^{-c} \cos(\log(e^{-\pi k}x^{\pi T})) -
\mbox{lin}_{c/T}(e^{-k}x^T)}{\log(e^{-k}x^T)}\,.$$
The complicated cosine term disappears at $x=e^{1/(2T)}$. Then
(\ref{d_thm1p2}) becomes
\begin{eqnarray}\label{d_thm1p3}
\Theta_cf(e^{1/(2T)})\,=\,\frac{4T}{\pi} \sum_{k\in\mathbb{Z}}
\frac{(-1)^{k+1}}{(2k-1)^2}\, e^{(k-1/2)c/T} f\left(e^{k/T}\right).
\end{eqnarray}
In order to obtain the $\Theta_c$ derivative of $f$ at $x$, we consider
the function
$g\,:\, t\mapsto f(xe^{-1/(2T)}t).$
It satisfies the assumptions used for deducing (\ref{d_thm1p3}). Now, applying
(\ref{d_thm1p3}) to $g$, we arrive at the desired formula (\ref{d_thm1.1}).

We still have to get rid of the additional assumption (\ref{d_thm1p1}).
If $f$ is any function in $B_{c,\pi T}^\infty$, then
$$ f_\varepsilon(x)\,:=\,f(x^{1-\varepsilon}) \,x^{c\varepsilon(T-1)} \mbox{lin}_c(x^{\varepsilon T})$$
belongs to $B_{c,\pi T}^\infty$ for each $\varepsilon \in (0, 1)$ and it satisfies
(\ref{d_thm1p1}). Applying (\ref{d_thm1.1}) to $f_\varepsilon$ and letting $\varepsilon \rightarrow 0_+$,
we find that (\ref{d_thm1.1}) holds for $f$ as well.
\hfill$\Box$
\vskip0,3cm

We note that formula (\ref{d_thm1.1}) yields a very short proof for a
Bernstein-type inequality for Mellin derivatives in $L^p$ norms for any
$p\in [1, \infty].$ Indeed, by the triangular inequality for norms we have
\begin{eqnarray}\label{Bernst1}
\left\|\Theta_cf\right\|_{X_c^p}\,\leq\, \frac{4T}{\pi} \sum_{k\in\mathbb{Z}}
\frac{1}{(2k-1)^2}\,e^{(k-1/2)c/T}\,\left\|f( \cdot \,e^{(k-1/2)/T})\right\|_{X_c^p}.
\end{eqnarray}
It is easily verified that for any positive $a$, there holds
$$\left\|f( \cdot\, a)\right\|_{X_c^p}\,=\, a^{-c} \|f\|_{X_c^p}\,.$$
Furthermore, it is known that
$$\sum_{k\in\mathbb{Z}} \frac{1}{(2k-1)^2}\,=\, \frac{\pi^2}{4}.$$
Thus, it follows from (\ref{Bernst1}) that
\begin{eqnarray}\label{Bernst2}
\|\Theta_cf\|_{X_c^p}\,\leq\, \pi T \|f\|_{X_c^p}\,.
\end{eqnarray}

Inequality (\ref{Bernst2}) in conjunction with Theorem \ref{derivative} and the Mellin
inversion formula shows that if $f\in B_{c,\pi T}^p$ for some $p\in[1, \infty]$,
then $\Theta_cf \in B_{c,\pi T}^p$ as well.

If $f$ does not belong to $B_{c,\pi T}^\infty$ but the two sides of
formula (\ref{d_thm1.1}) exist, we may say that (\ref{d_thm1.1}) holds with
a remainder $(R^B_{\pi T} f)(x)$ defined as the deviation of the right-hand
side from $\Theta_cf(x)$. We expect that $|(R^B_{\pi T}f)(x)|$ is small if
$\Theta_cf$ is close to $B_{c,\pi T}^\infty.$

For the Mellin inversion class for $p=2$,  $M_c^2$, we may state a
precise result as follows.

\begin{Theorem}\label{d_thm2}
Let $f\in \mathcal{M}_c^2$ and suppose that $v [f]_{M_c}^\wedge(c+iv)$ is absolutely
integrable on $\mathbb{R}$ with respect to $v$.
Then, for any $T>0$ and $x\in\mathbb{R}^+$, we have
\begin{eqnarray}\label{d_thm2.1}
\Theta_cf(x)\,=\, \frac{4T}{\pi} \sum_{k\in\mathbb{Z}} \frac{(-1)^{k+1}}{(2k-1)^2}
\,e^{(k-1/2)c/T} f\left(x e^{(k-1/2)/T}\right) + (R^B_{\pi T}f)(x),
\end{eqnarray}
where
\begin{eqnarray}\label{d_thm2.2}
(R^B_{\pi T}f)(x)\,=\, \frac{1}{2\pi i} \int_{|v|\geq\pi T}
\left[v- \pi T\phi\left(\frac{v}{\pi T}\right)\right] [f]_{M_c}^\wedge(c+iv) x^{-c-iv}dv
\end{eqnarray}
with
\begin{eqnarray}\label{d_thm2.3}
\phi(v)\,:=\, \bigg|v+1 -4\left\lfloor\frac{v+3}{4}\right\rfloor \bigg|-1.
\end{eqnarray}
In particular,
\begin{eqnarray}\label{d_thm2.4}
|(R^B_{\pi T}f)(x)|\,\leq\, \frac{x^{-c}}{2\pi} \int_{|v|\geq \pi T} 
\left(|v|+\pi T\right)|[f]_{M_c}^\wedge(c+iv)|dv
\end{eqnarray}
and
\begin{eqnarray}\label{d_thm2.5}
\|R^B_{\pi T}f\|_{X_c^\infty}\,\leq\, \frac{1}{\pi} 
\mbox{\rm dist}_1\left(\Theta_cf, B_{c,\pi T}^2\right).
\end{eqnarray}
\end{Theorem}
\begin{figure}[h]
\begin{center}
\setlength{\unitlength}{.6mm}
\begin{picture}(160,32)(-80,-15)
\linethickness{.4pt}
\put(-75,0){\vector(1,0){155}}
\put(0,-15){\vector(0,1){35}}
\multiput(-60,0)(20,0){7}{\line(0,-1){2}}
\thicklines
\multiput(-50,-10)(40,0){3}{\line(1,1){20}}
\multiput(-70,10)(40,0){4}{\line(1,-1){20}}
\put(77,2){\makebox(0,0)[b]{$v$}}
\put(2,20){\makebox(0,0)[tl]{$\phi(v)$}}
\put(-20,-3){\makebox(0,0)[t]{$-2\hphantom{-}$}}
\put(20,-3){\makebox(0,0)[t]{$2$}}
\end{picture}
\end{center}
\caption{\label{phi}
The graph of the function $\phi$.
}
\end{figure}
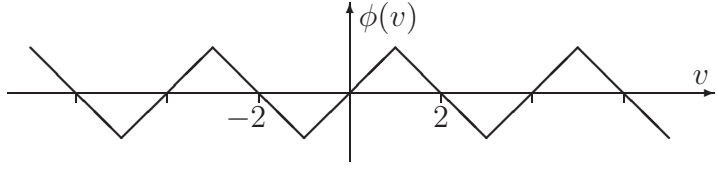
{\bf Proof.}\,
Define
\begin{eqnarray}\label{d_thm2p1}
f_1(x)\,:=\, \frac{1}{2\pi} \int_{|v|\geq \pi T} [f]_{M_c}^\wedge(c+iv)
x^{-c-iv}dv.
\end{eqnarray}
Then $f-f_1\in B_{c,\pi T}^\infty$, and so (\ref{d_thm1.1}) applies. It
yields that
\begin{eqnarray}\label{d_thm2p2}
(R^B_{\pi T}f)(x)\,=\, \Theta_c f_1(x) - \frac{4T}{\pi} \sum_{k\in\mathbb{Z}}
\frac{(-1)^{k+1}}{(2k-1)^2}\, e^{(k-1/2)c/T}\,f_1\left(xe^{(k-1/2)/T}\right).
\end{eqnarray}
We know from Theorem~\ref{derivative} that
$$\Theta_c f_1(x)\,=\, \frac{1}{2\pi i} \int_{|v|\geq \pi T} v
[f]_{M_c}^\wedge(c+iv) x^{-c-iv} dv.$$
Furthermore, by (\ref{d_thm2p1}),
$$f_1\left(x e^{(k-1/2)/T}\right)\,=\, \frac{1}{2\pi} \int_{|v|\geq \pi T}
[f]_{M_c}^\wedge(c+iv)\left( x e^{(k-1/2)/T}\right)^{-c-iv} dv.$$
Using these integral representations and interchanging summation and
integration, which is allowed by Levi's theorem, we may rewrite
(\ref{d_thm2p2}) as
\begin{eqnarray}\label{d_thm2p3}
(R^B_{\pi T} f)(x)\,=\, \frac{1}{2\pi i} \int_{|v|\geq \pi T}
\bigg(v -\psi(v)\bigg) [f]_{M_c}^\wedge(c+iv) x^{-c-iv} dv,
\end{eqnarray}
where
$$ \psi(v)\,:=\, \frac{4Ti}{\pi} \sum_{k\in\mathbb{Z}}
\frac{(-1)^{k+1}}{(2k-1)^2}\,
e^{-i(k-1/2)v/T}.$$
Now, for $v\in\mathbb{R}$, consider the function $g_v\,:\, x \mapsto ix^{-iv}.$
We note that $g_v \in B_{0,\pi T}^\infty$ if $|v|\leq \pi T.$
Hence $g_v$ satisfies the hypotheses of Theorem~\ref{d_thm1} for $c=0$ and
this restriction on $v$. Since $\Theta_0 g_v(1)=v$, we find by applying
(\ref{d_thm1.1}) to $g_v$ with $c=0$ and $x=1$ that $\psi(v)= v$ for
$v\in [-\pi T, \pi T].$

We also note that $\psi(v+ 2\pi T)=-\psi(v)$ and (consequently) $\psi(v+4\pi
T) = \psi(v).$ Hence $\psi$ is a $4\pi T$-periodic function that is given
on the interval $[-\pi t, 3\pi T]$ by
\begin{eqnarray*}
\psi(v)\, =\left\{
\begin{array}{cl}
v & \hbox{ if } -\pi T\le v \le \pi T,\\ \\
2\pi T -v & \hbox{ if } \pi T \leq v \leq 3\pi T.
\end{array}
\right.
\end{eqnarray*}
Thus, using the function $\phi$ defined in (\ref{d_thm2.3}),
 whose graph is shown in Fig.~\ref{phi},
 we can express $\psi(v)$ as $\pi T\phi(v/(\pi
T)).$ Hence (\ref{d_thm2p3}) implies (\ref{d_thm2.2}).

Inequalities (\ref{d_thm2.4}) and (\ref{d_thm2.5}) are easily obtained by noting that
$|\phi(v)|\leq 1$ for $v\in\mathbb{R}$ and by recalling Corollary~\ref{cor1}.
\hfill $\Box$

\subsection{An extension of the Bernstein-type inequality}

Just the same way as we deduced (\ref{Bernst2}) from (\ref{d_thm1.1}), we may use
(\ref{d_thm2.1}) to obtain
$$ \|\Theta_c f\|_{X_c^p}\,\leq\, \pi T \|f\|_{X_c^p} + \|R^B_{\pi T}
f\|_{X_c^p}.$$

For $p=2$ we can profit from the isometry of the Mellin transform
expressed by the formula
$$ \|f\|_{X_c^2}\,=\, \frac{1}{\sqrt{2\pi}}\left(\int_\mathbb{R} |[f]_{M_c}^\wedge(c+iv)|^2dv \right)^{1/2},$$
obtaining the following theorem
\begin{Theorem}\label{Bernapprox}
Under the assumptions of Theorem \ref{d_thm2} we have
$$\|\Theta_c  f\|_{X_c^2} \,\leq\,\pi T \|f\|_{X_c^2} + \frac{1}{\sqrt{2\pi}}
\,\mbox{\rm dist}_2(\Theta_c f, B_{c,\pi T}^2)$$
for any $T>0$.
\end{Theorem}
{\bf Proof}. 
With $f_1$ defined in (\ref{d_thm2p1}) and $f_0:=f-f_1$, we have
\begin{eqnarray}\label{Bernst3}
\|\Theta_c f\|_{X_c^2}\,\leq\, \|\Theta_c f_0\|_{X_c^2} +
\|\Theta_c f_1\|_{X_c^2} \,\leq\,
\pi T \|f_0\|_{X_c^2} +\|\Theta_c f_1\|_{X_c^2}
\end{eqnarray}
since (\ref{Bernst2}) applies to $f_0$.
We are going to estimate the quantities on the right-hand side in terms of $
f$.
Using the isometry of the Mellin transform, we find that
\begin{align*}
\|f\|_{X_c^2}^2 &=\, \frac{1}{2\pi} \int_\mathbb{R}|[f]_{M_c}^\wedge(c+iv)|^2 dv\\
&=\, \frac{1}{2\pi}\left[
\int_{|v|\leq \pi T}|[f]_{M_c}^\wedge(c+iv)|^2 dv +
\int_{|v|\geq \pi T}|[f]_{M_c}^\wedge(c+iv)|^2 dv\right] \\
&=\, \|f_0\|_{X_c^2}^2 +\|f_1\|_{X_c^2}^2\,,
\end{align*}  
which implies that $\|f_0\|_{X_c^2} \leq\|f\|_{X_c^2}.$ Next we note that
\begin{align*}
\|\Theta_c f_1\|_{X_c^2} &=\,
\frac{1}{\sqrt{2\pi}} \left(\int_\mathbb{R} |\left[\Theta_cf_1\right]_{M_c}^\wedge(c+iv)|^2dv\right)^{1/2}\\
&=\,\frac{1}{\sqrt{2\pi}} \left(\int_\mathbb{R} |v\left[f_1\right]_{M_c}^\wedge(c+iv)|^2
dv\right)^{1/2}\\
&=\,\frac{1}{\sqrt{2\pi}} \left(\int_{|v|\geq\pi T}|v[f]_{M_c}^\wedge(c+iv)|^2
dv\right)^{1/2}\\
&=\, \frac{1}{\sqrt{2\pi}} \mbox{dist}_2(\Theta_c f, B_{c,\pi T}^2).
\end{align*}
Thus (\ref{Bernst3}) implies the assertion. \hfill $\Box$
\vskip0,4cm
\noindent
{\bf Aknowledgments}. Carlo Bardaro and Ilaria Mantellini have been partially supported by the ``Gruppo Nazionale per l'Analisi Matematica e Applicazioni (GNAMPA) of the ``Istituto Nazionale di Alta Matematica'' (INDAM) as well as by the Department of Mathematics and Computer Sciences of the University of Perugia.



\begin{thebibliography}{1}
\bibitem{AD} R. A. Adams, Sobolev Spaces, Academic Press, New York, London, 1975.
\bibitem{BBM0} C. Bardaro, P.L. Butzer and I. Mantellini, The exponential sampling theorem of signal analysis and the reproducing kernel formula in the 
Mellin transform setting, Sampling Theory in Signal and Image Processing, 13(1), (2014), 35--66.
\bibitem{BBM} C. Bardaro, P.L. Butzer and I. Mantellini, The foundations of fractional calculus in the Mellin transform setting with applications, 
J. Fourier analysis and applications, 21 (2015), 961--1017.
\bibitem{BBM2} C. Bardaro, P.L. Butzer and I. Mantellini, The Mellin-Parseval formula and its interconnections with the exponential sampling theorem of 
optical physics, Integral Transforms and special functions, 27(1), (2016), 17--29.
\bibitem{BBMS} C. Bardaro, P.L. Butzer, I. Mantellini and  G. Schmeisser, On the Paley-Wiener theorem in the Mellin transform setting, Journal of Approximation Theory, to appear (2016).
\bibitem{BOA} R.P. Boas, The derivative of a trigonometric integral,
J. London Math. Soc. 12 (1937), 164--165.
\bibitem{BP} M. Bertero and E.R. Pike, Exponential sampling method for Laplace and other dilationally invariant transforms I. Singular-system analysis. II. 
Examples in photon correction spectroscopy and Fraunhofer diffraction, {Inverse Problems}, 7 (1991), 1--20, 21--41.
\bibitem{BJ0} P.L. Butzer and S. Jansche, Mellin transform theory and the role of its differential and integral operators, 
In ``Proc. Workshop on Transform Methods and Special Functions'', Varna, 1996, Inst. Math. and Inf., Bulgarian Acad. Sci., Sofia, 1998. 
\bibitem{BJ1} P.L. Butzer and S. Jansche, A direct approach to Mellin transforms, J. Fourier analysis and application, 3 (1997), 325--375.
\bibitem{BJ3} P.L. Butzer and S. Jansche, The exponential sampling theorem of signal analysis, Atti Sem. Mat. Fis. Univ. Modena, 
Suppl. Vol.\ 46,  99--122, (1998),  special issue dedicated to Professor Calogero Vinti.
\bibitem{BJ2} P.L. Butzer and  S. Jansche, A self contained approach to Mellin transform analysis for square integrable functions and applications, 
Integral Transforms and Special Functions, 8 (1999), 175--198.
\bibitem{BKT} P.L. Butzer,  A.A. Kilbas and  J.J. Trujillo, Mellin transform analysis and integration by parts for Hadamard-type fractional integrals, 
J. Math. Anal. Appl., 270, (2002), 1--15.
\bibitem{BN} P.L. Butzer and R.J. Nessel, Fourier Analysis and Approximation Vol I, Academic Press, New York, 1971.
\bibitem{BSS1} P.L. Butzer, G. Schmeisser and R.L. Stens, Shannon's sampling theorem for bandlimited signals and their Hilbert transform, Boas-type formulae 
for higher order derivatives---The aliasing error involved by their extensions from bandlimited to non-bandlimited signal, Entropy, 14, (2012), 2192--2226.
\bibitem{BSS2} P.L. Butzer, G. Schmeisser and R.L. Stens, Basic relations valid for the Bernstein Space $B^p_\sigma$ and their extensions to functions 
from larger spaces with error estimates in terms of their distances from $B^p_\sigma,$ J. Fourier Analysis and Applications, 19 (2013), 333--375.
\bibitem{BSS3} P.L. Butzer, G. Schmeisser and R.L. Stens, Basic relations valid for the Bernstein spaces $B^2_\sigma$ and their extensions to larger 
functions spaces via a unified distance concept, In  ``Function Spaces X'', Banach Centre Publs, Vol 102, Warszawa 2014, 41--55.
\bibitem{MA} R.G. Mamedov, The Mellin Transform and Approximation Theory, (in Russian), ``Elm'' Baku, 1991.
\bibitem{OSP} N. Ostrowsky, D. Sornette, P. Parker and E.R. Pike, Exponential sampling method for light scattering polydispersity analysis, 
{Opt. Acta}, {28} (1994),  1059--1070.
\bibitem{RIE} M. Riesz, Formule d'interpolation pour la d\'eriv\'ee d'un polyn\^ome,
Comptes Rendus Acad. Sci. Paris 158 (1914), 1152--1154.
\bibitem{SCH2} G. Schmeisser, Numerical differentiation inspired by a formula of R.\,P.\
Boas, J.\ Approx.\ Theory 160 (2009), 202--222.
\end{thebibliography}
\end{document}